\documentclass[12pt]{amsart}
\usepackage{graphicx, overpic, epsf,amssymb,amscd,times}

\newtheorem{thm}{Theorem}[section]

\newtheorem{prop}[thm]{Proposition}
\newtheorem{defn}[thm]{Definition}
\newtheorem{lemma}[thm]{Lemma}
\newtheorem{cor}[thm]{Corollary}

\newtheorem{q}[thm]{Question}

\newcommand{\R}{\mathbb{R}}
\newcommand{\Z}{\mathbb{Z}}
\newcommand{\Q}{\mathbb{Q}}

\newcommand{\bdry}{\partial}

\newcommand{\s}{\vskip.1in}
\newcommand{\n}{\noindent}

\newcommand{\be}{\begin{enumerate}}
\newcommand{\ee}{\end{enumerate}}
\newcommand{\veer}{Veer(S,\bdry S)} 
\newcommand{\dehn}{Dehn^+(S,\bdry S)}
\newcommand{\dehnnb}{Dehn^+(S)}

\begin{document}


\title{Right-veering diffeomorphisms of compact surfaces with boundary I}

\author{Ko Honda}
\address{University of Southern California, Los Angeles, CA 90089} \email{khonda@math.usc.edu}
\urladdr{http://rcf.usc.edu/\char126 khonda} 

\author{William H. Kazez}
\address{University of Georgia, Athens, GA 30602} \email{will@math.uga.edu}
\urladdr{http://www.math.uga.edu/\char126 will} 

\author{Gordana Mati\'c}
\address{University of Georgia, Athens, GA 30602} \email{gordana@math.uga.edu}
\urladdr{http://www.math.uga.edu/\char126 gordana}

\date{This version: October 27, 2005. (The pictures are in color.)}

 \thanks{KH supported by an Alfred P.\ Sloan Fellowship and an NSF CAREER Award; GM supported by NSF grant DMS-0072853; WHK supported by NSF grant DMS-0073029.} 

\begin{abstract}
We initiate the study of the monoid of right-veering diffeomorphisms on a compact oriented surface with nonempty boundary. The monoid strictly contains the monoid of products of positive Dehn twists. We explain the relationship to tight contact structures and open book decompositions.
\end{abstract}

\maketitle

\section{Introduction}

Let $S$ be a compact oriented surface with nonempty boundary. Denote by $Aut(S,\bdry S)$ the group of (isotopy classes of) diffeomorphisms of $S$ which restrict to the identity on $\bdry S$. (Such diffeomorphisms are automatically orientation-preserving.) In this paper we introduce the monoid $\veer\subset Aut(S,\bdry S)$ of {\em right-veering} diffeomorphisms of $S$, and explore how it is different from the monoid $\dehn\subset Aut(S,\bdry S)$ of products of positive Dehn twists. Informally said, a diffeomorphism $h\in Aut(S,\bdry S)$ is {\em right-veering} if every properly embedded arc $\alpha$ on $S$ is mapped ``to the right'' under $h$ -- see Section~\ref{defns} for a precise definition.

Our primary motivation for studying right-veering diffeomorphisms on $S$ comes from Giroux's work \cite{Gi} which relates contact structures and open book decompositions. Giroux showed that there is a 1-1 correspondence between isotopy classes of contact structures on a closed oriented 3-manifold $M$ and equivalence classes of open book decompositions on $M$ modulo a certain stabilization operation. We will denote an open book decomposition by $(S,h)$, where $S$ is a compact oriented surface with boundary (the {\em page}) and $h\in Aut(S,\bdry S)$ is the {\em monodromy map}. Our first theorem is the following:

\begin{thm} \label{veer}
A contact structure $(M,\xi)$ is tight if and only if all of its open book decompositions $(S,h)$ have right-veering $h\in Aut(S,\bdry S)$.
\end{thm}

\n
As usual, we assume that 3-manifolds are oriented, and contact structures are cooriented. Theorem~\ref{veer} is an improvement of the ``sobering arc'' criterion for overtwistedness given in Goodman's thesis \cite{Go}.

Theorem~\ref{veer} indicates that there is a gulf between $\veer$ and $\dehn$, in view of the characterization by Giroux of Stein fillable contact structures as those admitting $(S,h)$ with $h\in \dehn$.  In \cite{HKM2} we will show how to detect certain differences between $\veer$ and $\dehn$ without resorting to contact topology. 


\section{Right-veering} \label{defns}

Let $S$ be a compact oriented surface with nonempty boundary, and let $\alpha$ and $\beta$ be isotopy classes (relative to the endpoints) of properly embedded oriented arcs $[0,1]\rightarrow S$ with a common initial point $\alpha(0)=\beta(0)=x\in \partial S$. We will often blur the distinction between isotopy classes of arcs/curves and the individual arcs/curves if there is no danger of confusion.  Let $\pi:\tilde S\rightarrow S$ be the universal cover of $S$ (the interior of $\tilde S$ will always be $\R^2$ since $S$ has at least one boundary component), and let $\tilde x\in \bdry \tilde S$ be a lift of $x\in \partial S$. Take lifts $\tilde \alpha$ and $\tilde \beta$ of $\alpha$ and $\beta$ with $\tilde \alpha(0) =\tilde \beta(0)=\tilde x$. $\tilde \alpha$ divides $\tilde S$ into two regions -- the region ``to the left'' (where the boundary orientation induced from the region coincides with the orientation on $\tilde \alpha$) and the region ``to the right''. We say that $\beta$ is {\em to the right} of $\alpha$, denoted $\alpha \ge \beta$, if either $\alpha=\beta$ (and hence $\tilde\alpha(1)=\tilde\beta(1)$), or $\tilde\beta(1)$ is in the region to the right.

Alternatively, isotop $\alpha$ and $\beta$, while fixing their endpoints, so that they intersect transversely (this include the endpoints) and with the fewest possible number of intersections (we refer to this as intersecting {\em efficiently}). Assume that $\alpha\not=\beta$.  Then in the universal cover $\tilde S$, $\tilde\alpha$ and $\tilde\beta$ will meet only at $\tilde x$. If not, subarcs of $\tilde \alpha$ and $\tilde\beta$  would cobound a disk $D$ in $\tilde S$, and we could use an innermost disk argument on $\pi(D)\subset S$ to reduce the number of intersections of $\alpha$ and $\beta$ by isotopy. Then $\alpha\ge \beta$ if $int(\tilde\beta)$ lies in the region to the right. As an alternative to passing to the universal cover, we simply check to see if the tangent vectors $(\dot \beta(0), \dot\alpha(0))$ define the orientation on $S$ at $x$. 

\begin{defn}
Let $h:S \to S$ be a diffeomorphism that restricts to the identity map on $\partial S$. Let $\alpha$ be a properly embedded oriented arc starting at a basepoint $x \in \partial S$. Then $h$ is {\em right-veering} if for every choice of basepoint $x \in \partial S$ and every choice of $\alpha$ based at $x$, $h(\alpha)$ is to the right of $\alpha$ (at $x$). If $C$ is a boundary component of $S$, we say is $h$ is {\em right-veering with respect to $C$} if $h(\alpha)$ is to the right of $\alpha$ for all $\alpha$ starting at a point on $C$.
\end{defn}

\n
{\bf Remark.}
 We will tacitly assume that our arcs $\alpha$ are not boundary-parallel, although the inclusion of such arcs will not affect the definition of a right-veering diffeomorphism.

\s\n
{\bf Remark.} A notion equivalent to right-veering (see Theorem~\ref{monotone}) was considered by Amoros-Bogomolov-Katzarkov-Pantev in \cite{ABKP}. The notion of ``veering to the right'' in the universal cover can be traced at least as far back as Thurston's proof of the left orderability of the braid group (a fact originally due to Dehornoy~\cite{De}). Right-veering diffeomorphisms are called diffeomorphisms with ``protected boundary'' in Goodman \cite{Go}. In particular, Goodman shows that there are overtwisted contact structures which have monodromy maps which are right-veering; this is similar to our Proposition~\ref{rvstabilization}.

\s
We denote the set of isotopy classes of right-veering diffeomorphisms by $\veer\subset Aut(S,\bdry S)$. Since the composition of two right-veering diffeomorphisms is again right-veering, $\veer$ is a submonoid of $Aut(S,\bdry S)$.

Let us now interpret the notion of right-veering in terms of the circle at infinity. Assume the Euler characteristic $\chi(S)$ is negative, i.e., $S$ is not a disk or an annulus. We endow $S$ with a hyperbolic metric for which $\bdry S$ is geodesic. The universal cover $\pi: \tilde S\rightarrow S$ can then be viewed as a subset of the Poincar\'e disk $D^2 = \mathbb{H}^2 \cup S^1_\infty$. Now let $C$ be a component of $\bdry S$ and $L$ be a component of $\pi^{-1}(C)$. If $h\in Aut(S,\bdry S)$, let $\tilde h$ be the lift of $h$ that is the identity on $L$. Now the closure of $\tilde S$ in $D^2$ is a starlike disk. One portion of $\bdry \tilde S$ is $L$, and its complement in $\bdry \tilde S$ will be denoted $L_\infty$. Orient $L_\infty$ using the boundary orientation of $\tilde S$ and then linearly order the interval $L_\infty$ via an orientation-preserving homeomorphism with $\R$. The lift $\tilde h$ induces a homeomorphism $h_\infty:L_\infty \to L_\infty$. Also, given two elements $a,b$ in $Homeo^+(\R)$, the group of orientation-preserving homeomorphisms of $\R$, we write $a\geq b$ if $a(z)\geq b(z)$ for all $z\in \R$ and $a>b$ if $a(z)>b(z)$ for all $z\in \R$.  $a$ is said to be {\em nonincreasing} if $id\geq a$, and {\em strictly decreasing} if $id> a$. 

\begin{figure}[ht]
\begin{overpic}[height=4cm]{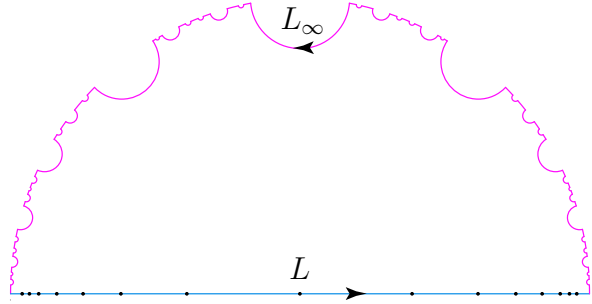}
\put(46.5,47){$L_{\infty}$}
\put(48,4){$L$}
\end{overpic}
\caption{$\tilde S$ as a subset of the Poincar\'e disk.}
\label{L_infinity}
\end{figure}

We have the following theorem, whose proof will occupy the next section: 

\begin{thm} \label{monotone}
Let $S$ be a hyperbolic surface with geodesic boundary and $h\in Aut(S,\bdry S)$. Then the following are equivalent: \be
\item $h$ is right-veering with respect to the boundary component $C$. \item $h$ sends every (properly) immersed geodesic arc $\alpha$ which begins on $C$ to the right.
\item $id\geq h_\infty$. (Here $h_\infty$ is defined with respect to $C$.) \ee
\end{thm}

Observe that our definition of $\alpha\geq \beta$ also makes sense for immersed geodesic arcs $\alpha$ and $\beta$.

The following is a useful way to show that certain diffeomorphisms are right-veering: 

\begin{lemma} \label{st}
Let $h\in Aut(S,\bdry S)$ be a right-veering diffeomorphism. Suppose $S'$ is obtained from $S$ by attaching a 1-handle $B$, and $h'=h\cup id_B$, namely $h'$ is the extension of $h$ by the identity map on $B$. Then $h'$ is right-veering.
\end{lemma}

\begin{proof}
Assume that $S$ is hyperbolic with geodesic boundary. (This excludes the case where $S$ is an annulus, but it is easy to furnish a proof in this case.)  Attaching a 1-handle $B$ is equivalent to gluing a pair-of-pants $P$ (taken to be hyperbolic with geodesic boundary) to $S$ along some common boundary components (1 or 2).  Let $S'=S\cup P$.  If $\alpha$ is a properly embedded essential arc of $S'$ with both endpoints on $\bdry S$, write it as a decomposition into subarcs (in order) $\alpha_1, \beta_1, \alpha_2,\beta_2,\dots,\beta_{n-1},\alpha_n$, where $\alpha_i$, $i=1,\dots,n$, are arcs on $S$ and $\beta_i$, $i=1,\dots,n-1$, are arcs in $P$. The case when $\alpha$ has at least one endpoint on $P$ is similar and will be omitted.  Since $h$ is right-veering on $S$, $h(\alpha_i)$ are to the right of $\alpha_i$.  We assume that $\alpha$ is piecewise geodesic, namely each $\alpha_i$ and $\beta_i$ is geodesic.

Now lift $\alpha$ and $h'(\alpha)$ to the universal cover $\pi:\tilde S'\rightarrow S'$ so that they have the same initial point $\tilde x\in \pi^{-1}(\bdry S)$.  Assume inductively that $\alpha_2*\beta_2*\dots*\beta_{n-1}*\alpha_n$ is to the left of $h'(\alpha_2)*h'(\beta_2)*\dots*h'(\beta_{n-1})*h'(\alpha_n)$.  Here $*$ denotes concatenation of arcs.  Then we prove that $\alpha$ is to the left of $h'(\alpha)$. If the lifts $\tilde\alpha_1$ and $\tilde h(\tilde \alpha_1)$ (with the same initial points) have terminal points on different connected components of $\pi^{-1}(\bdry S)$, then $\alpha$ is immediately to the left of $h'(\alpha)$.  If the terminal points lie on the same connected component of $\pi^{-1}(\bdry S)$, then either $\alpha_1=h(\alpha_1)$ (and the terminal points coincide) or the terminal point $\tilde y_1$ of $\tilde \alpha_1$ is to the left of the terminal point $\tilde y_2$ of $\tilde h(\tilde\alpha_1)$.  Now, using the inductive hypothesis, the lift of $\alpha_2*\beta_2*\dots*\beta_{n-1}*\alpha_n$ starting at $\tilde y_1$ is to the left of the lift of $h'(\alpha_2)*h'(\beta_2)*\dots*h'(\beta_{n-1})*h'(\alpha_n)$ starting at $\tilde y_1$, which in turn is to the left of the lift of $h'(\alpha_2)*h'(\beta_2)*\dots*h'(\beta_{n-1})*h'(\alpha_n)$ starting at $\tilde y_2$. 
\end{proof}

\begin{cor}
Let $S$ be a compact surface with nonempty boundary, obtained by gluing two compact surfaces $S_1$ and $S_2$ along some collection of boundary components in such a way that $\bdry S_i\cap \bdry S\not=\emptyset$ for $i=1,2$. If $h_i\in Aut^{\geq 0}(S_i,\bdry S_i)$, $i=1,2$, then $h= h_1\cup h_2 \in \veer$.
\end{cor}

\begin{proof}
The condition in the Corollary assures us that $S$ can be obtained from each of $S_i$, $i=1,2$, by consecutively attaching 1-handles. By Lemma~\ref{st}, $id_{S_1}\cup h_2$ and $h_1\cup id_{S_2}$ are both right-veering. Their composition is clearly right-veering as well.
\end{proof}

There is another monoid $\dehn$ which is well-known in symplectic geometry, namely the monoid of products of positive Dehn twists. More precisely, if $\gamma$ is a closed, homotopically nontrivial curve on $S$, let $R_\gamma$ denote  the positive (= right-handed) Dehn twist about $\gamma$. Then $h\in \dehn$ if and only if $h$ can be written as $R_{\gamma_1}\dots R_{\gamma_k}$, where $k$ may be zero. Although every element of the mapping class group on a closed surface can be written as a product of positive Dehn twists (hence $\dehnnb=Aut(S)$ if $S$ is closed), the same does not hold for $\dehn$ and $Aut(S,\bdry S)$ due to the existence of the boundary. We have the following: 

\begin{lemma}
$\dehn\subset \veer$.
\end{lemma}

\begin{proof}
Let $A$ be an annulus and $R_\gamma$ be a positive Dehn twist along the core of this annulus. It is easy to verify that $R_\gamma \in Aut(A,\bdry A)$ is right-veering. Now, by repeated application of Lemma~\ref{st}, we see that a positive Dehn twist $R_\gamma$ along an essential curve is right-veering, provided $\gamma$ is either (i) nonseparating or (ii) separates $S$ into two components, each of which nontrivially contains a component of $\bdry S$. Now, a positive Dehn twist $R_\gamma$ along a curve $\gamma$ that does not satisfy either (i) or (ii) can be written as a product of positive Dehn twists about nonseparating curves. Finally, observe that $\veer$ is closed under composition.
\end{proof}

One of the goals of this paper is to understand to what extent $\veer$ and $\dehn$ differ.

\section{Proof of Theorem~\ref{monotone}} 

In this section we give a proof of Theorem~\ref{monotone}. (In fact, we will prove much more!) Along the way we will discuss Thurston's classification of surface diffeomorphisms \cite{Th} (also see \cite{FLP, CB, Bo}), and then give a criterion for determining precisely when pseudo-Anosov homeomorphisms and periodic homeomorphisms are right-veering (Propositions~\ref{pseudo} and \ref{periodic}). Let us first dispense with the easy observations: 

(2)$\Rightarrow$(1) is immediate.

(3)$\Rightarrow$(2):
Suppose $z\ge h_\infty(z)$ for all $z\in L_\infty$. Given any point $x\in C$ and properly immersed geodesic arc $\alpha$ on $S$ with $\alpha(0)=x$, consider its lift $\tilde \alpha$ with $\tilde\alpha(0)=\tilde x$, where $\tilde x\in \tilde C_0 \cap \pi^{-1}(x)$. Then $\tilde h$ takes $\tilde\alpha(1)$ to $\widetilde{h(\alpha)}(1)$ (this follows from the definition of $\tilde h$). Now, $\tilde\alpha(1) \ge \tilde h (\tilde\alpha(1))= \widetilde{h(\alpha)}(1)$. Hence $\alpha\ge h(\alpha)$.

\subsection{Thurston classification of surface homeomorphisms} We now describe Thurston's classification, which improved upon earlier work of Nielsen \cite{Ni1, Ni2, Ni3}. A diffeomorphism $h:S\rightarrow S$ is {\em reducible} if there exists an {\em essential multicurve} $\gamma$, none of whose components is parallel to a component of $\bdry S$, such that $h(\gamma)=\gamma$. Here, a {\em multicurve} is a union of pairwise disjoint simple closed curves, and a multicurve is {\em essential} if no component of $S-\gamma$ is a disk or an annulus. If $h$ is not isotopic to a {\em reducible} diffeomorphism, then it is said to be {\em irreducible}, and is {\em freely} isotopic a homeomorphism $\psi$ of one of the following two types:

\be
\item A {\em periodic} homeomorphism, i.e., there is an integer $n>0$ such that $\psi^n=id$. \item A {\em pseudo-Anosov} homeomorphism. \ee

\subsection{Pseudo-Anosov case}
Suppose that $S$ is a hyperbolic surface with geodesic boundary (in other words, $S$ is not an annulus). We will first consider the pseudo-Anosov case. A pseudo-Anosov homeomorphism $\psi$ is equipped with the stable and unstable measured geodesic laminations $(\Lambda^s,\mu^s)$ and $(\Lambda^u,\mu^u)$ (here $\Lambda^s$ and $\Lambda^u$ are the laminations and $\mu^s$ and $\mu^u$ are the transverse measures) such that $\psi(\Lambda^s) =\Lambda^s$ and $\psi(\Lambda^u)=\Lambda^u$. $\Lambda$ (= $\Lambda^s$ or $\Lambda^u$) is {\em minimal} (i.e., does not contain any sublaminations), does not have closed or isolated leaves, is disjoint from the boundary $\bdry S$, and every component of $S-\lambda$ is either an open disk or a semi-open annulus containing a component of $\bdry S$. In particular, every leaf of $\Lambda$ is dense in $\Lambda$. We also have $\psi(\Lambda^s,\mu^s)=(\Lambda^s,\tau \mu^s)$ and $\psi(\Lambda^u,\mu^u)=(\Lambda^u,\tau^{-1}\mu^u)$ for some $\tau>1$. 

Let $C$ be a boundary component of $S$. Then the connected component of $S-\Lambda^s$ containing $C$ is a semi-open annulus $A$ whose metric completion $\hat{A}$ has geodesic boundary consisting of $n$ infinite geodesics $\lambda_1,\dots,\lambda_n$. Suppose that the $\lambda_i$ are numbered so that $i$ increases (modulo $n$) in the direction given by the boundary orientation on $C$. Now let $P_i\subset A$ be a semi-infinite geodesic which begins on $C$, is perpendicular to $C$, and runs parallel to $\gamma_i$ and $\gamma_{i+1}$ (modulo $n$) along the ``spike'' that is ``bounded'' by $\gamma_i$ and $\gamma_{i+1}$. These $P_i$ will be referred to as the {\em prongs}. Let $x_i=P_i\cap C$ be the endpoint of $P_i$ on $C$. We may assume that $\psi$ permutes (rotates) the prongs and, in particular, that there exists an integer $k$ so that $\psi: x_i\mapsto x_{i+k}$ for all $i$. 

Now, given a diffeomorphism $h\in Aut(S,\bdry S)$, let $H: S\times[0,1]\rightarrow S$ be an isotopy from $h$ to its pseudo-Anosov representative $\psi$. Define $\beta: C \times [0,1] \to C\times [0,1]$ by sending $(x,t)\mapsto (H(x,t),t)$. It follows that the arc $\beta(x_i \times [0,1])$ connects $(x_i, 0)$ and $(x_{i+k}, 1)$, where $k$ is as above. We call $\beta$ a {\it fractional Dehn twist} by an amount $c\in \Q$, where $c\equiv k/n$ modulo $1$ is the number of times $\beta(x_i \times [0,1])$ circles around $C \times [0,1]$ (here circling in the direction of $C$ is considered positive). Form the union of $C \times [0,1]$ and $S$ by gluing $C \times \{1\}$ and $C\subset \bdry S$. By identifying this union with $S$, we construct the homeomorphism $\beta\cup\psi$ on $S$ which is isotopic to $h$, relative to $\bdry S$. (We will assume that $h=\beta\cup\psi$, although $\psi$ is usually just a homeomorphism, not a diffeomorphism.)
For a surface with multiple boundary components $C_j$, we produce fractional Dehn twists $\beta_j$ by amounts $c_j$.

\s\n
{\bf Remark.} The difference between the automorphism $h\in Aut(S,\bdry S)$ and its psuedo-Anosov representative $\psi$ with respect to a boundary component was originally exploited by Gabai. (See \cite{GO} for the relationship to essential laminations and Dehn fillings and \cite{Ga} for a discussion and open questions degeneracy slope of a knot.)

\begin{prop} \label{pseudo}
If $h$ is isotopic to a pseudo-Anosov homeomorphism, then the following are equivalent:
\be
\item $h$ is right-veering with respect to $C$. \item $c>0$ for the boundary component $C$ . \item $id\geq h_\infty$ for the boundary component $C$. \ee
\end{prop}

\begin{proof}
We use the notation from the paragraph preceding the statement of Theorem~\ref{monotone}. Suppose $S$ is hyperbolic with geodesic boundary. Let $C$ be a component of $\bdry S$ for which the fractional Dehn twist is by an amount $c$. Let $P_1,\dots, P_n$ be the geodesic prongs that end on $C=C\times\{0\}$. Their union is preserved by $\psi$ (on its domain of definition), and the $P_i$ are subjected to the fractional Dehn twist $\beta$ along $C\times[0,1]$. For convenience, for any $j\in \Z$ we write $P_j$ to indicate $P_i$, where $i\equiv j$ modulo $n$ and $1\leq i\leq n$.

Suppose that $c>0$. Then $h_\infty:L_\infty\rightarrow L_\infty$ maps the terminal point of $\tilde P_i$ to the terminal point of $\tilde h (\tilde P_i)$, which is the same as the terminal point of $\tilde P_{i+cn}$.  (Here $\tilde P_j$, $j\in\Z$, is a lift of $P_j$ with initial point on $L$, where $j$ increases in the direction given by the orientation on $L$, which in turn is induced by the orientation on $\tilde S$.) If we denote the terminal point of $\tilde P_i$ in $L_\infty$ by $a_i$, then the interval $[a_i,a_{i+1})$ gets mapped to $[a_{i+cn},a_{i+cn+1})$ by the monotonicity. Therefore, $h_\infty$ is a strictly decreasing function (i.e., $id> h_\infty$), and all three notions of right-veering in Theorem~\ref{monotone} coincide in this case. 

Similarly, if $c<0$, then $id < h_\infty$ and $h$ moves every properly embedded essential arc $\alpha$ on $S$ strictly to its left. 

Now suppose that $c=0$. The pseudo-Anosov property implies that the $a_i$ are attracting fixed points for $\psi_\infty: L_\infty\rightarrow L_\infty$. It follows that there are points $z\in L_\infty$ for which $z\geq h_\infty(z)$ and others for which $z\leq h_\infty(z)$. To prove the proposition (i.e., to show $h$ is neither left- or right-veering), however, we need to do better -- we need to find a properly embedded arc $\alpha$ in $S$ which starts at $C$ and is sent to the left by $h$. We construct $\alpha$ as follows: Take an embedded arc $\alpha_1$ which starts at $C$ and transversely crosses an infinite boundary leaf of $\bdry \hat A$, say $\lambda_1$, and ends on a nearby leaf $\lambda$. (For example, it could be a subarc of a prong for $\Lambda^u$.) Let $\beta$ be the subarc of $\alpha_1$ from $\lambda_1$ to $\lambda$. At $\lambda \cap \alpha_1$, veer left, and follow $\lambda$ until the first return to $\beta$ -- call this arc $\alpha_2$. The existence of a point of first return is guaranteed by the minimality of $\Lambda^s$ ($\lambda$ is dense in $\Lambda^s$). At $\alpha_2\cap \beta$, either turn to the left or to the right so that we return to $C$ by moving parallel to $\alpha_1$ -- this gives us $\alpha_3$. Now let $\alpha=\alpha_1*\alpha_2*\alpha_3$.  We observe that $\alpha$ is embedded and that $\alpha$ is not boundary-parallel. Now, by the pseudo-Anosov property, $h$ contracts the interval $\beta$ by a factor strictly between 0 and 1, and moves $\lambda$ strictly to the left along $\beta$. Hence $h(\alpha)$ is either equal to $\alpha$ or $h(\alpha)\geq \alpha$. To eliminate the possibility $h(\alpha)=\alpha$, observe that a pseudo-Anosov map cannot fix any closed curves, in particular $\gamma$ obtained by closing up $\alpha_2$ with a subarc of $\beta$.
\end{proof}

\subsection{Periodic case}
Next consider the case when $h$ is freely isotopic to a periodic map $\psi$. As in the pseudo-Anosov case, we consider the trace $\beta: \bdry S\times[0,1] \rightarrow \bdry S\times[0,1]$ of an isotopy from $h$ to $\psi$ on $\partial S$ and write $h=\beta \cup \psi$. Let $C_0,\dots,C_k$ be the boundary components of $S$, and $c_i\in \Q$, $i=1,\dots, k$, be the amount of boundary twisting on the component $C_i$.

\begin{prop} \label{periodic}
If $h$ is freely isotopic to a periodic homeomorphism, then the following are equivalent:
\be
\item $h$ is right-veering with respect to the boundary component $C_0$.

\item $c_0> 0$ or else $c_0=0$ and $c_i\geq 0$ for all $i=1,\dots,k$.

\item $id\geq h_\infty$ for the boundary component $C_0$.
\ee
\end{prop}

\begin{proof}
We use the same notation as before. First observe that if $x<y$ on $L_\infty$, then $h_\infty(x)<h_\infty(y)$ since $h_\infty$ is orientation-preserving. Hence if $z>h_\infty(z)$ (resp.\ $z<h_\infty(z)$), then $z> h^n_\infty(z)$ (resp.\ $z<h^n_\infty(z)$). We therefore conclude that $z>h_\infty(z)$ (resp.\ $z=h_\infty(z)$, $z<h_\infty(z)$) is equivalent to $z>h^n_\infty(z)$ (resp.\ $z=h^n_\infty(z)$, $z<h^n_\infty(z)$) for any $n>0$. 

Next observe that, since $\psi$ is periodic, there exists a positive integer $n$ for which $\psi^n=id$ and $h^n= R_{\gamma_0}^{n_0}\dots R_{\gamma_k}^{n_k}$, where $\gamma_i$ is a parallel, disjoint copy of $C_i$, the $\gamma_i$ are pairwise disjoint, and $n_i\in \Z$. Moreover, $n_i= c_i n$. 

Suppose $c_0>0$, or, equivalently, $n_0>0$. We claim that $z>h_\infty(z)$ for all $z\in L_\infty$. In fact, take a properly embedded arc $\beta$ starting at $C=C_0$ and ending at $C_i$, $i\in\{0,\dots,k\}$, and consider its lifts $\tilde \beta_j$, $j\in \Z$, which start at $L$. Here $j$ increases in the direction given by the orientation of $L$. Now let $L_j$ be the component of $\pi^{-1}(C_i)$ which contains the endpoint $\tilde\beta_j(1)$. Observe that $L_j$ and $L_{j+1}$ are distinct.  (If not, then $L_i=L_j$ for all $i,j$, which leads to the contradiction that $L_i=L_j=L$.) Then $h^n_\infty$ maps $L_j$ to $L_{j+n_0}$ and moves the interval between $L_j$ and $L_{j+n_0}$ to the right of $L_{j+n_0}$. (In other words, $h^n_\infty$ is a strict monotonically decreasing function.) Hence $id > h^n_\infty$ and therefore $id>h_\infty$. 

Similarly, if $c_0<0$, then $id<h_\infty$. 

Next suppose that $c_0=0$ (or $n_0=0$). If $L'\subset L_\infty$ is a component of $\pi^{-1}(C_i)$ for some $i$, then $h^n_\infty$ maps $L'$ to itself. If all the $n_i\geq 0$, then $h^n_\infty$ either fixes the points on $L'$ or moves them to the right. Since the union of $(\cup_{i=0}^k\pi^{-1}(C_i))\cap L_\infty$ is dense in $L_\infty$, it follows that $id\geq h^n_\infty$ and $id\geq h_\infty$. On the other hand, if $n_i<0$ for some $i$, then we can easily find a properly embedded arc $\beta$ from $C=C_0$ to $C_i$ which is sent strictly to the left via $h^n$. If $h^n(\beta)> \beta$, then $h(\beta) > \beta$.
\end{proof}

\subsection{Reducible case}
Suppose $h\in Aut(S,\bdry S)$ is reducible. Then there exists a (nonempty) essential multicurve $\gamma$ which satisfies $h(\gamma)=\gamma$. Assume that $\gamma$ is maximal among all such multicurves. Fix a component $C$ of $\bdry S$. Then let $S_\gamma$ be the component of the metric closure of $S-\gamma$ which contains $C=C_0$. The other boundary components of $S_\gamma$ will be called $C_1,\dots, C_k$. $h$ may permute the components of $C_i$, $i=1,\dots,k$. Now, by the Nielsen-Thurston classification, $h|_{S_\gamma}$ is either periodic or pseudo-Anosov.

Suppose that $h$ is right-veering. We will show that $id\geq h_\infty$. But first we have the following lemma: 

\begin{lemma}\label{subsurface}
Let $S$ be a hyperbolic surface with geodesic boundary and $g\in Aut(S,\bdry S)$. Let $S'\subsetneq S$ be a subsurface, also with geodesic boundary, and let $C$ be a common boundary component of $S$ and $S'$. If $g$ is the identity map when restricted to $S'$, $\delta$ is a closed curve parallel to and disjoint from $C$, and $m$ is a positive integer, then $id> (R^m_\delta g)_\infty$ with respect to $C$.
\end{lemma}

\begin{proof}
Let $L$ be the component of $\pi^{-1}(C)$ which is fixed under the lift $\tilde g$. Let $\alpha$ be a properly embedded essential arc in $S'$ starting on $C$ and ending on $C'\not = C$. By hypothesis, $\tilde g$ is also the identity on the lifts $\tilde\alpha_j$, $j\in \Z$, of $\alpha$ with initial points on $L$. Order $\tilde\alpha_j$ along $L$ so that $j$ increases in the direction given by the orientation of $L$. Let $L_j$ the component of $\pi^{-1}(C')$ that $\tilde\alpha_j$ ends on. Then $\widetilde{R_\delta^m g}$ maps $L_j$ to $L_{j+m}$. (Note that the $L_j$ are distinct.) This implies that any point on $L_\infty$ between $L_j$ and $L_{j+m}$ gets mapped to the right of $L_{j+m}$. Hence $id> (R_\delta^m g)_\infty$.
\end{proof}

We state the following corollary of Lemma~\ref{subsurface}, which will become useful in Section~\ref{stab}:

\begin{cor}\label{rvtest}
Let $S$ be a hyperbolic surface with geodesic boundary and $g\in Aut(S,\bdry S)$. Let $C$ be a component of $\partial S$ and $C'$ be an embedded closed geodesic in $S$ (which also may be a component of $\bdry S$). Suppose we have the following: (i) $g(C')=C'$, and (ii) there is a properly embedded arc $\alpha$ that starts on $C$ and ends on $C'$ so that $g(\alpha)$ is isotopic to $\alpha$ (via an isotopy which takes $C'$ to itself but does not necessarily fix each point of $C'$). If $\delta$ is a closed curve parallel to and disjoint from $C$, then $id> (R_\delta g)_\infty$ with respect to $C$.
\end{cor}

\begin{proof}
Consider the union $C\cup\alpha\cup C'$. Its neighborhood $S'$ is a pair-of-pants on which we may take $h$ to be the identity after isotopy (relative to $\bdry S$). In particular, any Dehn twists about curves parallel to $C'$ can be pushed away from $S'$. Now apply Lemma~\ref{subsurface}.
\end{proof}

Let us consider the case $h|_{S_\gamma}= R_\delta^m$, where $\delta$ is a closed curve parallel to and disjoint from $C$, and $m\in \Z$. Suppose first that $m=0$. By appealing to the universal cover, it is not difficult to see that $id\geq h_\infty$ if and only if $id \geq (h|_{S-S_\gamma})_\infty$ with respect to any of the components $C_i$, $i=1,\dots,k$. It follows that if $h$ is right-veering with respect to $C$, then $h|_{S-S_\gamma}$ is also right-veering with respect to $C_i$, $i=1,\dots, k$. We may inductively excise such $S_\gamma$ with $h|_{S_\gamma}=id$ from $S$. Next suppose that $m>0$. Then we are in the case of Lemma~\ref{subsurface} and $id> h_\infty$. If $m<0$, then $id< h_\infty$. 

Now assume that $h|_{S_\gamma}$ is not the identity and is either periodic or pseudo-Anosov. Note that by the maximality of $\gamma$, there are no arcs from $C_0$ to $C_i$ which are preserved by $h$.  Let $L$ be a lift of $C_0$ in the universal cover $\tilde S_\gamma$, and $L_\infty$ be the complement of $L$ in $\bdry \tilde S_\gamma$. If $\alpha$ is an arc from $C_0$ to $C_i$, then the terminal point of its lift $\tilde\alpha$ must lie on a different component of $\pi^{-1}(\cup C_i)$ from that of $\tilde h(\tilde\alpha)$. Since $h(\alpha)$ cannot be to the left of $\alpha$ (this would contradict the right-veering property for any extension of $\alpha$ to $S$), $h|_{S_\gamma}$ is right-veering.  It remains to show that $h|_{S_\gamma}$ is right-veering implies $id\ge h_\infty$.

Suppose $h|_{S_\gamma}$ is pseudo-Anosov. Then there exist prongs $P_i$ in $S_\gamma$ which end along $C$ as in the proof of Proposition~\ref{pseudo}, and $h$ rotates the prongs by $c>0$, according to Proposition~\ref{pseudo}. This fact immediately forces $h_\infty$ to be strictly decreasing. Note that Proposition~\ref{pseudo} is still valid even when only $C=C_0$ is fixed and the $C_i$, $i=1,\dots,k$, are permuted. 

Next suppose $h|_{S_\gamma}$ is periodic and not the identity. Recalling that $id\geq h^n_\infty$ if and only if $id\geq h_\infty$, we only deal with $h^n$. If the amount of rotation is $c>0$, then $h^n_\infty$ is strictly decreasing. Suppose now that $c=0$. Then we may take $h^n$ to be the identity on $S_\gamma$ after isotopy. Again, we can excise such $S_\gamma$, thereby inductively shrinking $S$. This completes the proof of Theorem~\ref{monotone}.

\section{Open book decompositions}

In the section we review the fundamental work of Giroux (building on work of Thurston-Winkeln\-kemper \cite{TW}, Bennequin \cite{Be}, Eliashberg-Gromov \cite{EG}, and Torisu \cite{To}), which relates contact structures and open book decompositions. 

Let $(S,h)$ be a pair consisting of an oriented surface $S$ and a diffeomorphism $h:S \rightarrow S$ which restricts to the identity on $\partial S$, and let $K$ be a link in a closed oriented 3-manifold $M$. An {\em open book decomposition} for $M$ with {\em binding} $K$ is a homeomorphism between $((S \times[0,1]) /_{\sim _h}, (\partial S \times[0,1]) /_{\sim_h})$ and $(M,K)$. The equivalence relation $\sim_h$ is generated by $(x,1) \sim_h (h(x),0)$ for $x\in S$ and $(y,t) \sim_h (y,t')$ for $y \in \partial S$. We will often identify $M$ with $(S \times[0,1]) / _{\sim _h}$; with this identification $S_t= S \times \{t\}, t\in [0,1]$, is called a {\em page} of the open book decomposition and $h$ is called the {\em monodromy map}. Two open book decompositions are {\em equivalent} if there is an ambient isotopy taking binding to binding and pages to pages. We will denote an open book decomposition by $(S,h)$, although, strictly speaking, an open book decomposition is determined by the triple $(S,h,K)$. There is a slight difference -- if we do not specify $K\subset M$, we are referring to isomorphism classes of open books instead of isotopy classes. 

Every closed 3-manifold has an open book decomposition, but it is not unique. One way of obtaining inequivalent open book decompositions is to perform a {\em positive} or {\em negative stabilization}: $(S',h')$ is a {\em stabilization} of $(S,h)$ if $S'$ is the union of the surface $S$ and a band $B$ attached along the boundary of $S$ (i.e., $S'$ is obtained from $S$ by attaching a 1-handle along $\bdry S$), and $h'$ is defined as follows. Let $\gamma$ be a simple closed curve in $S'$ ``dual'' to the cocore of $B$ (i.e., $\gamma$ intersects the cocore of $B$ at exactly one point) and let $id_B \cup h$ be the extension of $h$ by the identity map to $B \cup S$. Then for a {\em positive} stabilization $h'=R_\gamma \circ (id_B \cup h)$, and for a {\em negative} stabilization $h'=R_\gamma^{-1} \circ (id_B \cup h)$. It is well-known that if $(S',h')$ is a positive (negative) stabilization of $(S,h)$, and $(S,h)$ is an open book decomposition of $(M,K)$, then $(S',h')$ is an open book decomposition of $(M, K')$ where $K'$ is obtained by a Murasugi sum of $K$ (also called the {\em plumbing} of $K$) with a positive (negative) Hopf link. 

A contact structure $\xi$ is said to be {\em supported} by the open book decomposition $(S,h,K)$ if there is a contact 1-form $\alpha$ which: \be
\item induces a symplectic form $d\alpha$ on each fiber $S_t$; \item $K$ is transverse to $\xi$, and the orientation on $K$ given by $\alpha$ is the same as the boundary orientation induced from $S$ coming from the symplectic structure.
\ee
Thurston and Winkelnkemper \cite{TW} showed that any open book decomposition $(S,h,K)$ of $M$ supports a contact structure $\xi$. Moreover, the contact planes can be made arbitrarily close to the tangent planes of the pages (away from the binding).

The following result is the converse (and more), due to Giroux~\cite{Gi}. 

\begin{thm}[Giroux]
Every contact structure $(M,\xi)$ on a closed 3-manifold $M$ is supported by some open book decomposition $(S,h,K)$. Moreover, two open book decompositions $(S,h,K)$ and $(S',h',K')$ which support the same contact structure $(M,\xi)$ become equivalent after applying a sequence of positive stabilizations to each.
\end{thm}

In the framework of open book decompositions, the {\em holomorphically fillable} contact structures admit a simple characterization, as a result of work by Loi-Piergallini \cite{LP}, combined with the relative version of the above theorem. (Also see \cite{AO}.) 

\begin{cor}[Loi-Piergallini, Giroux] \label{holo} A contact structure $\xi$ on $M$ is holomorphically fillable if and only if  $\xi$ is supported by some open book $(S,h,K)$ with $h\in \dehn$.
\end{cor}

\n
{\bf Remark.} Note that the result does not state that {\em any} open book $(S,h,K)$ for a holomorphically fillable contact structure has monodromy in $\dehn$, just that there is at least one such. It is an open question to determine whether all supporting open books $(S,h,K)$ for a holomorphically fillable $(M,\xi)$ have monodromy in $\dehn$.

\s\n
Our Theorem~\ref{veer} is a characterization of tightness in this framework of open book decompositions and contact structures. The ultimate goal is to be able to determine whether $(M,\xi)$ is tight, fillable, etc., just by looking at a single $(S,h)$ -- from that perspective, we still fall short of the goal....

\section{Proof of Theorem~\ref{veer}}

In this section we prove Theorem~\ref{veer}. The proof is based on the following idea: Recall that a bypass (cf. \cite{H1}) can be thought of as half of an overtwisted disk. Suppose we glue two contact manifolds $M_1$ and $M_2$ along a common boundary $\Sigma$. The most elementary way for
$M_1\cup M_2$ to be overtwisted is for a Legendrian arc in $\Sigma$ to be an arc of attachment for a bypass on each side of $\Sigma$. This might appear to be a very special kind of obstruction to a tight gluing of the two contact manifolds, and that we can match bypasses from $M_1$ and $M_2$ only in very special cases. However, there is a simple operation called {\em bypass rotation} which enables us to match up bypasses from $M_1$ and $M_2$ if $h:S\rightarrow S$ is not right-veering.

\begin{figure}[ht]
\begin{overpic}[height=3cm]{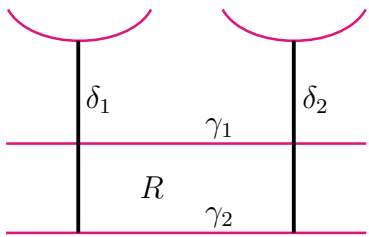}
\put(22,35){$\delta_1$}
\put(82,35){$\delta_2$}
\put(55,28.2){$\gamma_1$}
\put(55,3.7){$\gamma_2$}
\put(37,10){$R$}
\end{overpic}
\caption{$\delta_2$ is to the right of $\delta_1$, viewed from the interior of $N$.} \label{rotating}
\end{figure}

We now discuss {\em bypass rotation}, which first appeared in \cite{HKM3}. Let $\Sigma$ be a convex boundary component of a contact 3-manifold $(N,\zeta)$.  Although slightly awkward, we give $\Sigma$ the opposite orientation to the induced orientation from $N$. (In particular, in Figure~\ref{rotating}, the boundary orientation points into the page.)  Let $\delta_1$ and $\delta_2$ be disjoint Legendrian arcs in $\Sigma$ that have three transverse intersections with the dividing set $\Gamma_\Sigma$, including their endpoints. The arcs $\delta_1$ and $\delta_2$ are arcs of attachment for potential bypasses, which are to be attached ``from the front''. Suppose there is an embedded rectangle $R$, where two of the sides are subarcs of $\delta_1$ and $\delta_2$, and the other two sides are subarcs $\gamma_1$ and $\gamma_2$ of $\Gamma_\Sigma$. Assume $\delta_1$ and $\delta_2$ both start on $\gamma_2$, extend beyond $\gamma_1$ and do not reintersect $\bdry R$. If $R$, $\delta_1$, and $\delta_2$ are as in Figure~\ref{rotating}, that is, with the orientation induced from $R$ (which in turn is induced from that of $\Sigma$), $\gamma_1$ starts on $\delta_2$ and ends on $\delta_1$, then we say that $\delta_1$ lies {\it to the left} of $\delta_2$.   Note that Figure~\ref{rotating} only represents the local picture near $R\cup \delta_1\cup \delta_2$; in particular, any of the four dividing arcs may be part of the same dividing curve in $\Gamma_\Sigma$.

\begin{lemma}[Bypass rotation]\label{rightbypassrotating} Let $(N,\zeta)$ be a contact
3-manifold with convex boundary, and let $\delta_1$ and $\delta_2$ be disjoint arcs on a boundary component $\Sigma$ of $N$. If $\delta_1$ is to the left of $\delta_2$, and $\delta_2$ is an arc of attachment for a bypass inside $N$, then $\delta_1$ is also an arc of attachment for a bypass inside $N$.
\end{lemma}

\begin{proof}
Suppose there exists a bypass $B\subset N$ along $\delta_2$. Then attaching $B$ onto $\Sigma$ modifies $\Gamma_\Sigma$ to $\Gamma_{\Sigma'}$ so that $\delta_1$ becomes a trivial arc of attachment. By the Right-to-Life Principle \cite{H2,HKM} (i.e., bypasses always exist along trivial arcs of attachment), there exists a bypass along $\delta_1$.
\end{proof}

The ability to rotate bypasses will be coupled with the following lemma: 

\begin{lemma}\label{generalrotating}
If $\alpha$ and $\beta$ are oriented, properly embedded arcs in $S$ such that $\alpha \ge \beta$, then there exists a sequence of oriented, properly embedded arcs $\alpha = \alpha_0
\ge \dots \ge \alpha_n = \beta$ where $\alpha_i$ and $\alpha_{i+1}$ have disjoint interiors for $i =0, \dots, n-1$. Here we require the initial points $\alpha_i(0)$ to be the same, but the terminal points $\alpha_i(1)$ do not need to coincide.
\end{lemma}

\begin{proof}
Let us denote by $\#(\alpha,\beta)$ the geometric intersection number of $\alpha$ and $\beta$, with the exception of the endpoints. (Hence $\#(\alpha,\beta)=0$ if they have no interior intersection points.) We will show that, given $\alpha$ and $\beta$ in $S$ for which $\alpha \ge \beta$, $\alpha\not=\beta$ and $\#(\alpha,\beta)\not =0$, there always exists an arc $\alpha'$ such that $\alpha \ge \alpha'\ge \beta$ and $\#(\alpha,\alpha')<\#(\alpha,\beta)$, $\#(\alpha',\beta) < \#(\alpha,\beta)$. By recursively applying the above procedure, we construct the sequence $\alpha = \alpha_0 \ge \dots \ge \alpha_n = \beta$ with $\#(\alpha_i,\alpha_{i+1})=0$ for $i=1,\dots,n-1$. 

Suppose $\alpha \ge \beta$ and $\alpha\not=\beta$. We may assume that $\alpha$ and $\beta$ intersect transversely and efficiently. Then let $\alpha \cap \beta =\{p_0, p_1, \dots , p_n\} = \{q_0, q_1, \dots , q_n\}$, where $p_i=\alpha(t_i), 0=t_0 \le t_1 \le \dots \le t_n < 1$ and $q_i=\beta(s_i), 0=s_0 \le s_1 \le \dots \le s_n<1$. We will analyze several cases depending on the nature of the intersection of $\alpha$ and $\beta$ at $p_1$: 

\s\n
{\bf Case 1.} Assume that $p_1= q_r$ and at $p_1$ the tangent vectors to $\beta$ and $\alpha$ (in that order) determine the orientation on $S$. (See Figure~\ref{Cases1and2}).

Let $\alpha'$ be the curve obtained by following $\alpha$ until $p_1$ and then veering right and following $\beta$ from that point on; in other words, $\alpha'=\alpha|_{[0,t_1]}*\beta|_{[s_r,1]}$. Then $\#(\alpha,\alpha')=\#(\alpha,\beta)-r$, and $\#(\alpha',\beta)=0$. 

Since it is an important point, we will explain why $\alpha \ge \alpha'\ge \beta$. $\alpha'\ge \beta$ because $\alpha'$ is to the left of $\beta$ near their common initial point and $\alpha'$ and $\beta$ never reintersect (hence $\alpha'$ and $\beta$ intersect efficiently). $\alpha\ge \alpha'$ since $\alpha'$ starts to the right of $\alpha$ and smoothing the piecewise geodesic arc $\alpha|_{[0,t_1]}*\beta|_{[s_r,1]}$ into a (smooth) geodesic arc with the same endpoints pushes $\alpha'$ further to the right.

\begin{figure}[ht]
\begin{overpic}[height=5cm]{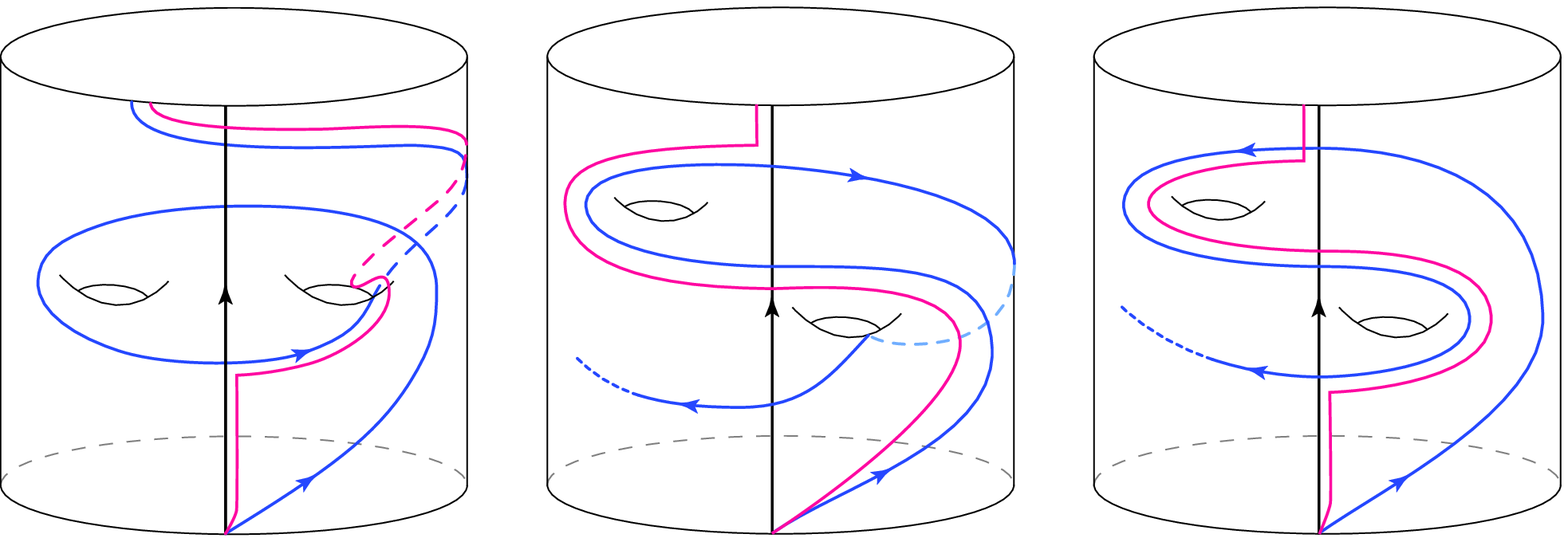}
\put(10,0){Case~1}
\put(11,11){$\alpha$}
\put(25,12){$\beta$}
\put(11,20){$p_1$}
\put(16,20){$q_r$}
\put(44,0){Case~2A}
\put(46,11){$\alpha$}
\put(61.5,12){$\beta$}
\put(45,28){$p_{r'}$}
\put(51,28){$q_{r''}$}
\put(79,0){Case~2B}
\put(81,11){$\alpha$}
\put(95,12){$\beta$}
\put(80.5,28.5){$p_{r'}$}
\put(86,28.5){$q_{r''}$}
\end{overpic}
\caption{}
\label{Cases1and2}
\end{figure}

\s\n
{\bf Case 2.} Assume that $p_1= q_r$, $r>1$, and at $p_1$ the local intersection is such that the tangent vectors to $\alpha$ and $\beta$ determine the orientation on $S$.

Let $p_{r'}$ be the last point on $\alpha$ where $\alpha$ intersects $\beta|_{[0,s_r]}$.  Say  $p_{r'}=q_{r''}$.  Case~2 is split into two subcases depending on the orientation of the curves at  $p_{r'}=q_{r''}$.

\s\n
{\bf Case 2A.} First suppose the tangent vectors to $\beta$ and $\alpha$ determine the orientation on $S$ at $p_{r'}=q_{r''}$.  Then let $\alpha'= \beta|_{[0,s_{r''}]} * \alpha|_{[t_{r'},1]}$, that is, we follow $\beta$ until we hit $p_{r'}$ and veer left and travel along $\alpha$ for the rest of the journey. (See Figure~\ref{Cases1and2}). It is clear that $\#(\alpha,\alpha')<\#(\alpha,\beta)$ and $\#(\alpha',\beta)< \#(\alpha,\beta)$. We claim that $\alpha\geq \alpha'\geq \beta$. Assuming $\alpha$ and $\beta$ are geodesics, we pass to the universal cover $\pi:\tilde S\rightarrow S$. Let $\tilde \alpha$ and $\tilde\beta$ be lifts of $\alpha$ and $\beta$ which start at the same basepoint $\tilde x\in \bdry \tilde S$. Then the lift $\tilde \alpha'$ of $\alpha'$ which starts at $\tilde x$ is a piecewise geodesic arc which follows $\tilde\beta$ until it hits a lift of $p_{r'}$ and goes left along a lift of $\alpha$. $\tilde \alpha'$ clearly ends to the left of $\tilde\beta$ and hence $\alpha'\geq \beta$. Now, $\tilde\alpha'$ starts to the right of $\tilde\alpha$ (along $\tilde \beta$), and then switches to a different component of $\pi^{-1} (\alpha)$ from $\tilde \alpha$. Hence $\tilde \alpha'$ ends to the right of $\tilde \alpha$ and $\alpha\geq \alpha'$. (Whether the switch is to the right or to the left is not important -- rather, it's the fact that we switch to a different component of $\pi^{-1}(\alpha)$. Distinct lifts of $\alpha$ never self-intersect!)


\s\n
{\bf Case 2B.}
Next suppose that the tangent vectors to $\alpha$ and $\beta$ determine the orientation on $S$ at $p_{r'}=q_{r''}$.  Then let $\alpha'= \alpha|_{[0,t_1]}* (\beta|_{[s_{r''}, s_r]})^{-1}* \alpha|_{[t_{r'},1]}$, i.e., we follow $\alpha$ until we hit $p_1$, veer right and travel backwards along $\beta$ until we reach $p_{r'}$, and then veer left and travel along $\alpha$ until the end. (See Figure~\ref{Cases1and2}).  Again, it is clear that $\#(\alpha,\alpha'), \#(\alpha',\beta) < \#(\alpha,\beta)$. To see that $\alpha\geq\alpha'\geq \beta$, consider the universal cover $\tilde S$ as in the previous paragraph. $\tilde\alpha'$ follows $\tilde \alpha$ until it hits a lift of $p_1$, goes right and follows a lift of $\beta$ until it hits another lift of $\alpha$ {\em besides} $\tilde \alpha$. Therefore, $\alpha\geq \alpha'$. Similar considerations give $\alpha'\geq \beta$.

\s\n
{\bf Case 3.} Assume that at $p_1= q_r$ the tangent vectors to $\alpha$ and $\beta$ determine the orientation on $S$ and that $r=1$. 

\s\n
{\bf Case 3A.} Assume that $\gamma=\alpha|_{[0,t_1]} \cup \beta|_{[0,s_1]}$ separates $S$. Let $R\subset S$ be the subsurface cut off by $\gamma$ (the one which does not intersect $\beta|_{(s_r,1)}$). 

Suppose first that there is a boundary component of $S$ in $R$. In this case, we take $\alpha'$ to be an arc in $R$ which connects $p_0=q_0$ to $\bdry S\cap R$ without intersecting the interior of $\alpha$ or $\beta$. Then $\#(\alpha,\alpha')=\#(\alpha',\beta)=0$ and it is evident that $\alpha\geq \alpha'\geq \beta$.

If there are no boundary components of $S$ in $R$, then there is an arc $\delta\subset R$ which starts at $p_0$, runs over a handle in $R$ and ends at $p_1$, and whose interior does not intersect $\alpha$ or $\beta$. Let $\alpha'= \delta* \beta|_{[s_1,1]}$. In this case, $\#(\alpha',\beta)=0$ and $\alpha'$ is clearly to the left of $\beta$. Although $\alpha'$ is to the right of $\alpha$, we apparently have not gained anything since we still have $\#(\alpha,\alpha')=\#(\alpha,\beta)$.  This is not a problem, as we will now be in Case~3B where $\gamma$ is nonseparating. See Figure~\ref{Cases3and4}. 

\begin{figure}[ht]
\begin{overpic}[height=5cm]{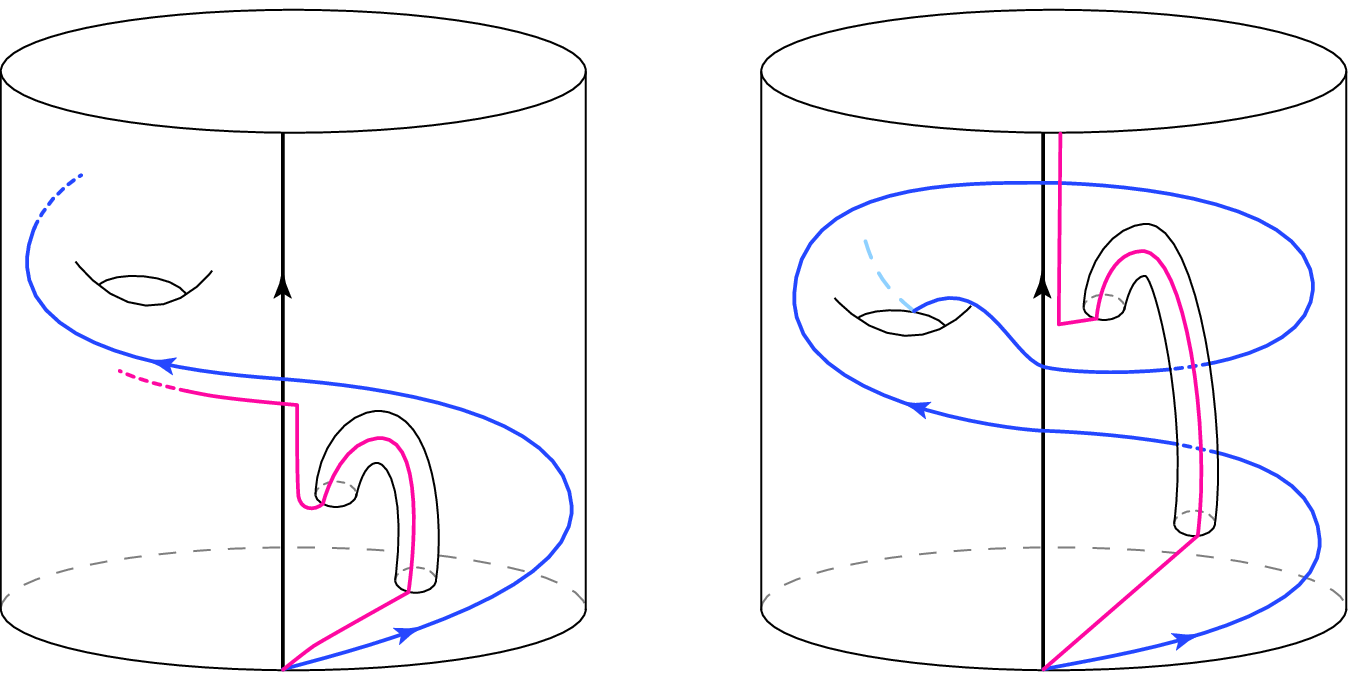}
\put(14,0){Case~3A}
\put(17,11){$\alpha$}
\put(23.5,11){$\delta$}
\put(39,11){$\beta$}
\put(16.5,30.5){$p_1$}
\put(23,30.5){$q_1$}
\put(71,0){Case~3B}
\put(73.5,11){$\alpha$}
\put(85.5,11){$\delta$}
\put(96.5,11){$\beta$}
\put(73,21.5){$p_1$}
\put(79.5,21.5){$q_1$}
\end{overpic}
\caption{}
\label{Cases3and4}
\end{figure}

\s\n
{\bf Case 3B.} Assume that $\gamma=\alpha|_{[0,t_1]}\cup \beta|_{[0,s_1]}$ is nonseparating.

Since $\gamma$ is nonseparating, there is an arc $\delta$ which connects $\alpha|_{[0,t_1]}$ to $\alpha|_{[t_i,t_{i+1}]}$ for some $i\ge 1$, and whose interior does not intersect $\beta$.
Construct an arc $\alpha'$ such that $\alpha \ge \alpha'\ge \beta$ by starting between $\alpha$ and $\beta$, following $\delta$ until $p_{i+1}=\alpha(t_{i+1})$ is reached, and then following $\alpha$ from there on. Then $\#(\alpha,\alpha')=0$, $\#(\alpha',\beta) \le \#(\alpha,\beta)-1$ and $\alpha\geq \alpha'\geq \beta$. See Figure~\ref{Cases3and4}.
\end{proof}

\begin{proof}[Proof of Theorem~\ref{veer}] An open book decomposition $(S,h,K)$ gives rise to a special Heegard decomposition of $M$. Denote by $H_i$, $i=1,2$, the two handlebodies $H_1=S \times[0,{1\over 2}] /_{\sim _h}$ and $H_2= S \times[{1\over 2}, 1] /_{\sim _h}$. Writing $S_t=S\times\{t\}$ as before, the oriented boundaries of $H_1$ and $H_2$ are $\Sigma_{1}=S_{1/2} \cup -S_0$ and $\Sigma_2=-S_{1/2} \cup S_{1}$, respectively.  See Figure~\ref{monodromy}. The orientation-reversing gluing map $g:\Sigma_2 \to \Sigma_1$ is given by $id_{S_{1/2}} \cup h$, where the monodromy map $h$ is regarded as a map of $S_1$ to $S_0$. Observe that, if the neighborhood $N(K)=K\times D^2$ of $K$ has coordinates $(z,r,\theta)$ and $S_t\cap (K\times D^2)$ is $\theta= 2\pi t$, then the separating surface $\Sigma\stackrel{def}=\Sigma_1=\Sigma_2$ of the Heegaard decomposition is smooth. 

\begin{figure}[ht]
\begin{overpic}[height=11cm]{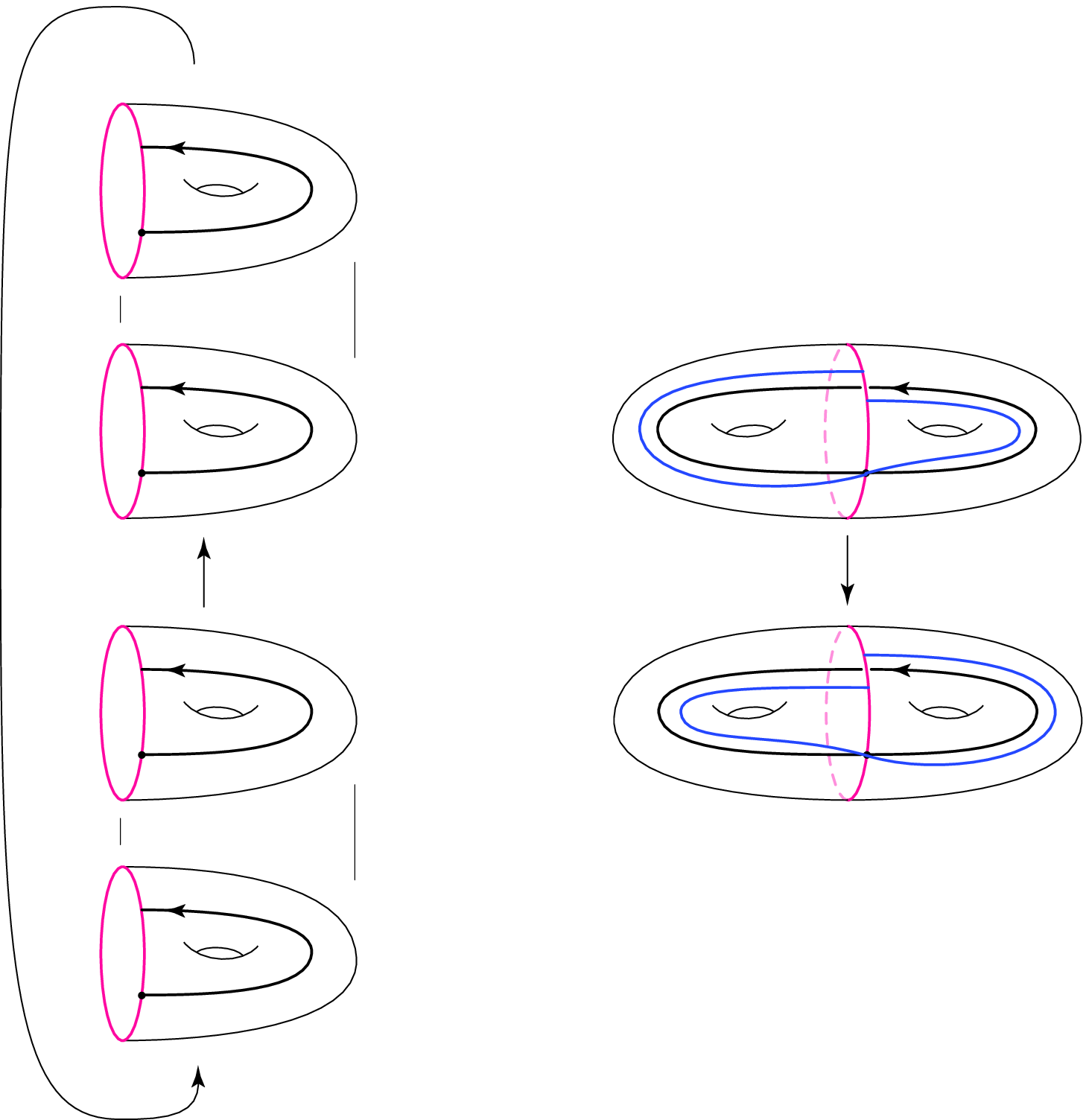}

\put(25,10.5){$\alpha$}
\put(4,24.3){$H_1$}
\put(4,71.5){$H_2$}
\put(1,48){$h$}
\put(20,48){$id_{S_{1/2}}$}
\put(33,13.5){$S_0$}
\put(33,35.3){$S_{1/2}$}
\put(33,60.7){$S_{1/2}$}
\put(33,82){$S_1$}

\put(49,35.3){$\Sigma_1$}
\put(49,60.7){$\Sigma_2$}

\put(63,25){$-S_0$}
\put(82,25){$S_{1/2}$}
\put(65,72){$S_1$}
\put(81,72){$-S_{1/2}$}

\put(71,42){$\gamma_1$}
\put(71,63){$\gamma_2$}
\put(71,35.7){$\delta_1$}
\put(71,68.3){$\delta_2$}

\put(68,48){$g=$}
\put(78,48){$h \cup id_{S_{1/2}}$}

\end{overpic}
\caption{}
\label{monodromy}
\end{figure}

Since $(S,h)$ is adapted to $\xi$, $\Sigma$ is a convex surface whose dividing set $\Gamma_\Sigma$ is equal to the binding $K$. Let $\alpha \subset S$ be an oriented, properly embedded arc and denote by $\alpha_t$ the copy of $\alpha$ in $S_t= S \times \{t\}$. The boundaries of the compressing disks $\alpha \times [0,1/2]\subset H_1$ and $\alpha \times [1/2,1] \subset H_2$ are closed curves $\gamma_1= \alpha^{-1}_0*\alpha_{1/2}$ and $\gamma_2=\alpha^{-1}_1*\alpha_{1/2}$, which we (not so) secretly view as parametrized by an interval. Now consider the arc $\delta_i$, obtained by ``perturbing the endpoints of $\gamma_i$'' as follows: Let $\alpha^L$ (resp.\ $\alpha^R$)
be an oriented, properly embedded arc on $S$ which has the same initial point as $\alpha$, is parallel to and disjoint from $\alpha$ with the exception of the initial point, and has final point slightly to the left (resp.\ right) of $\alpha$. Then let $\delta_1= (\alpha^L)^{-1}_0*\alpha_{1/2}$ and $\delta_2 =(\alpha^R)^{-1}_1* \alpha_{1/2}$. Apply the Legendrian realization principle, we may assume that $\delta_i\subset \Sigma$ is Legendrian. 

We now claim that $\delta_i$ is an arc of attachment for a bypass in $H_i$. The closed curve $\gamma_i$ intersects the dividing set $\Gamma_{\Sigma_i}$ in exactly two points corresponding to the two endpoints of $\alpha$. Take a parallel copy $\gamma^*_i$ of $\gamma_i$ (disjoint from $\delta_i$) on $\Sigma$. Apply the Legendrian realization principle to both $\gamma^*_i$ and $\delta_i$ (their union is still nonisolating), and split $H_i$ along a convex disk $D_i$ with boundary $\gamma^*_i$. (See Figure~\ref{decomp}.) Since $\#(\Gamma_{\Sigma_i}\cap\gamma^*_i)=2$, the dividing set $\Gamma_{D_i}$ is uniquely determined single arc. After rounding the edges, we observe that $\delta_i$ is a trivial arc of attachment on the cut-open manifold $H_i\setminus D_i$. The existence of a bypass along $\delta_i$ now follows from the Right-to-Life Principle. Alternatively we can argue that a convex handlebody $H$ which admits a collection of compressing disks $D_1,\dots,D_g$ satisfying $\#(\Gamma_{\bdry H}\cap \bdry D_i)=2$ and $H\setminus (D_1\cup \dots\cup D_g)=B^3$ is a standard neighborhood of a Legendrian graph, and the sought-after bypass is equivalent to the stabilization of a Legendrian arc. 

\begin{figure}[ht]
\begin{overpic}[height=4cm]{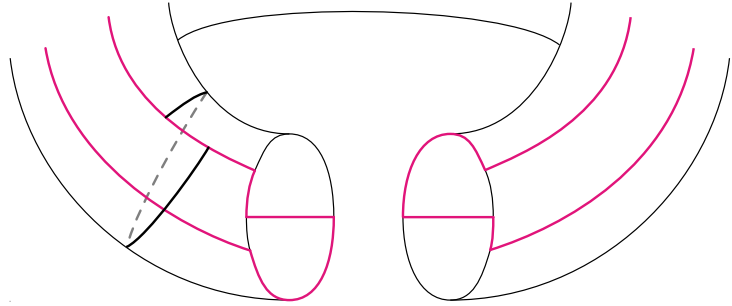}
\end{overpic}
\caption{The disk decomposition and the trivial bypass.} \label{decomp}
\end{figure}

Suppose $h$ is not right-veering. Then there exists an oriented, properly embedded arc $\alpha \subset S$ for which $h(\alpha)$ is to the left of $\alpha$. Suppose first that $h(\alpha)$ and $\alpha$ have disjoint interiors. We start with the two arcs of attachment $\delta_1$ and $\delta_2$, along which there are bypasses from $H_1$ and $H_2$. We view $\delta_2$ on $\bdry H_1$ by mapping via $g = h \cup id_{S_{1/2}}$ to $g(\delta_2) = (h(\alpha^R))^{-1}_0* \alpha_{1/2}$. Then $g(\delta_2)$ will be to the right of $\delta_1$ with respect to the boundary orientation on $\bdry H_1$. Since $g(\delta_2)$ is an arc of attachment for a bypass coming from $H_2$, it follows from Lemma~\ref{rightbypassrotating} that $\delta_1$ also is an arc of attachment for a bypass in $H_2$. The union of the two bypasses is an overtwisted disk, which contradicts the assumption that $(M,\xi)$ is tight. (Observe that if $h$ is the identity map, then $g(\delta_2)$ will be slightly to the left of $\delta_1$ with respect to the boundary orientation on $\bdry H_1$, and we are not able to construct an overtwisted disk by gluing the bypasses. This is good, since the contact manifold is the connected sum of $S^1\times S^2$'s with the unique tight contact structure!)

Now, if $h(\alpha)$ is to the left of $\alpha$ but the interiors are not disjoint, use Lemma~\ref{generalrotating} to produce a sequence $\alpha_0=h(\alpha) \ge \alpha_1 \ge \dots \ge \alpha_n=\alpha$ so that $\alpha_i$ and $\alpha_{i+1}$ have disjoint interiors for $i =0, \dots, n-1$. We use this sequence to construct a sequence of arcs of attachment $(\alpha_i)^{-1}_{0} * \alpha_{1/2}$ on $\partial H_1$ where consecutive arcs are disjoint, and $(\alpha_i)^{-1}_{0} * \alpha_{1/2}$ is to the right of $(\alpha_{i+1})^{-1}_{0} * \alpha_{1/2}$ (with respect to the boundary orientation on $\bdry H_1$). By applying a sequence of bypass rotations, it follows that $\delta_1$ bounds a bypass in $H_2$, giving rise to an overtwisted disk in $(M, \xi)$. 

Suppose now that $(M,\xi)$ is overtwisted. We will show that there is some open book decomposition $(S,h)$ for $(M,\xi)$ where $h$ is not right-veering. By Eliashberg's classification of overtwisted contact structures by their homotopy type (see \cite{El}), $(M,\xi)$ is isotopic to the connected sum $(M,\xi)\# (S^3,\xi^{OT})$, where $\xi^{OT}$ is the unique overtwisted contact structure on $S^3$ which is homotopic to the standard one, and the connected sum is done along a convex sphere $S^2$ with one dividing curve. Now, by \cite{DGS}, we may write $(S^3,\xi^{OT})$ as a connected sum $(S^3,\xi')\#(S^3,\xi'')$, where $(S^3,\xi')$ has an open book decomposition $(S',h')$ with $S'$ an annulus, and $h'$ is a negative Dehn twist about the core curve of the annulus. Performing a connected sum with $(S^3,\xi')$ is equivalent to performing a negative stabilization along a boundary-parallel (trivial) arc. Such an open book clearly does not have right-veering monodromy.
\end{proof}

\section{Stabilizing monodromy} \label{stab} 

In this section we prove the following proposition: 

\begin{prop}\label{rvstabilization}
Every $(S,h)$ can be made right-veering after a sequence of positive stabilizations.
\end{prop}

The idea behind the proof of Proposition~\ref{rvstabilization} is to create an example that is so strongly right-veering that when it is added to the boundary components of any $(S,h)$, the resulting surface and diffeomorphism are forced to be right-veering.

\begin{figure}[ht]
\begin{overpic}[ height=4.5cm]{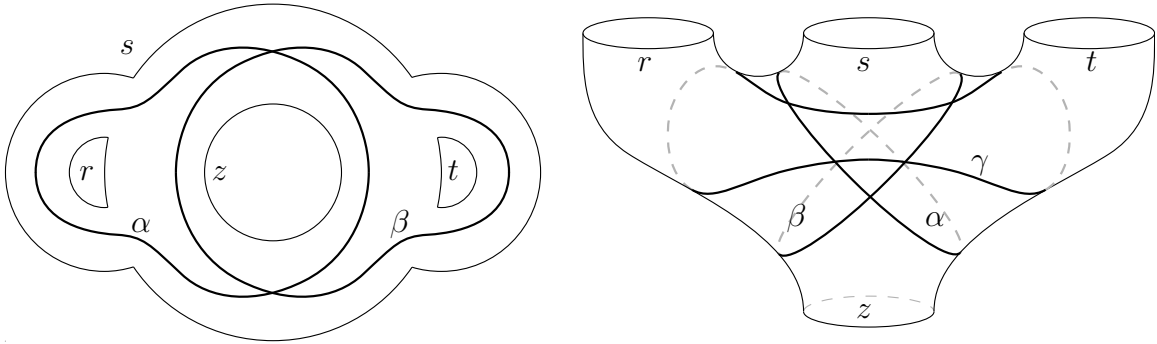}
\put(6.5,14){$r$}
\put(38.5,14){$t$}
\put(18,14){$z$}
\put(10,25){$s$}
\put(11,9.5){$\alpha$}
\put(33.5,9.5){$\beta$}
\put(55, 23.5){$r$}
\put(74, 23.5){$s$}
\put(94, 23.5){$t$}
\put(74, 2){$z$}
\put(80,10){$\alpha$}
\put(68,10){$\beta$}
\put(84,15){$\gamma$}
\end{overpic}
\caption{Two stabilizations (left) and the lantern relation (right)} \label{lantern}
\end{figure}

\s\n
{\bf Example.}
Let $A$ be an annulus with monodromy $id_A$. Denote the boundary components of $A$ by $x$ and $z$. Choose two properly embedded, boundary-parallel arcs $b_1$ and $b_2$ in $A$ such that $\bdry b_1, \bdry b_2\subset x$, and $b_1$ and $b_2$ intersect efficiently in two interior points. Thus there exist subarcs of $b_1$ and $b_2$ whose union is a core of $A$. Create $(L,h)$ by first stabilizing $(A, id_A)$ across $b_1$, and then stabilizing across $b_2$. (By ``stabilizing across an arc $b$'', we mean that the 1-handle $B=c\times [-\varepsilon,\varepsilon]$ is attached to the surface $S$ so that $\bdry b= \bdry c\times\{0\}$ and the extra positive Dehn twist is performed along the curve $b\cup (c\times\{0\})$.)
The result is that $L$ has three boundary components $r,s,t$ replacing $x$, and it still has the original boundary component $z$ of $A$. See Figure~\ref{lantern}.

By construction, $h$ is the composition of two positive Dehn twists $R_\alpha$ and $R_\beta$. Choose a curve $\gamma$ as shown in Figure~\ref{lantern}, and then by the lantern relation, $R_\gamma R_\beta R_\alpha = R_r R_s R_t R_z$. (Observe that $R_r$, $R_s$, $R_t$ and $R_z$ mutually commute.) It follows that $h= R^{-1}_\gamma R_r R_s R_t R_z$. 

We will now exhibit arcs $a_1$, $a_2$ and $a_3$ to which we apply Corollary~\ref{rvtest}.  Let $a_1$ be an embedded arc in $L$ that starts on $r$, ends on $t$, and does not intersect $\gamma$. From the lantern relation, we see that $a_1$ is preserved by $h$, modulo boundary twisting. (More precisely, $R^{-1}_r R^{-1}_t h$ fixes $a_1$, up to isotopy relative to $\partial L$.) Let $a_2$ be an embedded arc in $L$ that starts on $s$ and ends on $\gamma$. Again, the lantern relation shows that $R^{-1}_s R_\gamma h$ fixes $a_2$, up to isotopy relative to $s \cup \gamma$. Finally, an arc $a_3$ analogous to $a_2$ exists from $z$ to $\gamma$. This shows that $h$ is right-veering: $a_1$ implies that $h$ is right-veering with respect to $r$ and $t$, $a_2$ with respect to $s$, and $a_3$ with respect to $z$.

\s
The proof of Proposition~\ref{rvstabilization} is now immediate: 

\begin{proof}[Proof of Proposition~\ref{rvstabilization}] Let $(L, R_\beta R_\alpha)$ be as described in the above example. To each boundary component $z$ of $S$, attach a copy of $(L, R_\beta R_\alpha)$ along $z$. The surface produced by this procedure is evidently obtained by a sequence of positive stabilizations, and its monodromy is right-veering by Corollary~\ref{rvtest}.
\end{proof}

In a sense, Proposition~\ref{rvstabilization} is disappointing since any $(S,h)$ corresponding to an overtwisted contact structure can be disguised as a right-veering diffeomorphism by ``protecting the boundary''. 

The following questions arise naturally out of Proposition~\ref{rvstabilization}. 

\begin{q}
If $(S, h)$ is right-veering and cannot be destabilized, does it necessarily correspond to a tight contact structure? What if we assume that $h$ is pseudo-Anosov in addition (i.e., $h$ is isotopic to a pseudo-Anosov homeomorphism with only positive fractional Dehn twists)?
\end{q}

\begin{q} \label{pa}
If $(S,h)$ is right-veering and pseudo-Anosov, does it necessarily correspond to a tight contact structure?  
\end{q}

In other words, can you obtain a right-veering pseudo-Anosov homeomorphism by stabilizing an overtwisted $(S,h)$ as in Proposition~\ref{rvstabilization} or taking a further stabilization?

A related question is:

\begin{q}
Suppose $(S,h)$ is pseudo-Anosov, with fractional Dehn twists $c_i>0$. Determine an optimal constant $C$ (may be $\infty$) so that all the $(S,h)$ satisfying $c_i\geq C$ are tight.  Same for universally tight. Does $C$ depend on $S$? What are the conditions on $c_i$ (besides $c_i>0$) for the contact structure to be tight?
\end{q}

\begin{q}
Consider the contact manifold $((S\times[0,1])/\sim,\xi)$, where $(x,1)\sim (h(x),0)$ and $\xi$ is a tight contact 2-plane field which is close to the leaves. For which Dehn fillings is the resulting (canonical) contact manifold $(M,\xi)$ tight? 
\end{q}

Recall that open book decompositions are Dehn fillings of a special type.

We conclude this section by making some observations when $(S,h)$ is a stabilization.

\begin{lemma}
The open book decomposition $(S,h)$ is a (positive) stabilization if and only if there is a properly embedded arc $\alpha$ in $S$ such that $h(\alpha)$ intersects $\alpha$ only at the boundary and $h(\alpha)$ is to the right of $\alpha$ with respect to both endpoints of $\alpha$.
\end{lemma}

\begin{proof}
Assume $\alpha$ is as in the statement of the lemma. It is clear that $h(\alpha)= R_C(\alpha)$ if we denote by $C$ the class of the closed curve obtained as the union of $\alpha$ and $h(\alpha)$. Let us define $g=R_C^{-1} h$. Since $g(\alpha)=\alpha$, $g$ restricts to a homeomorphism that fixes the boundary of the surface $\Sigma=S \setminus \alpha$ obtained by cutting $S$ along $\alpha$. It is now easy to see that $(S,h)$ is a (positive) stabilization of $(\Sigma,g)$.
\end{proof}

\begin{lemma}
Let $S=T^2 \setminus D^2$ be a punctured torus and $(S,h)$ an open book decomposition with $h\in \veer$.  If $(S,h)$ is a stabilization, then $h\in \dehn$.
\end{lemma}

\begin{proof}
Since $(S,h)$ is a stabilization, there is an arc $\alpha$ such that $h(\alpha)$ intersects $\alpha$ only at its endpoints. Let $C$ be the class of the closed curve obtained as the union of $\alpha$ and $h(\alpha)$. Then $(S,h)$ is stabilization of $(A,g)$ for an annulus $A=S\setminus \alpha$. Then $g=R_{\gamma}^{ q}$, where $\gamma$ is the core of the annulus. It is easy to see that the composition $R_C  R_{\gamma}^{q}$ is not right-veering unless $q \geq 0$.  In fact, if $q<0$ and $\beta$ is an essential arc on $S$ which is parallel to $C$ (i.e., does not intersect $C$), then $h(\beta)> \beta$.
\end{proof}

\section{Contrasting right-veering and product of right Dehn twists}

By now the following theorem should be evident:

\begin{thm}   \label{notsame}
$\dehn \subsetneqq \veer$.
\end{thm}                                                                                       

In view of Theorem~\ref{monotone} (namely that $h\in Aut(S,\bdry S)$ right-veering is 
equivalent to $id\geq h_\infty$), Theorem~\ref{notsame} gives a negative answer to a
conjecture of Amor\'os-Bogomolov-Katzarkov-Pantev \cite{ABKP}.  (They had
conjectured that $\dehn= \veer$.)

The following is a list of right-veering diffeomorphisms $h\in \veer$
which are not in $\dehn$.

\begin{enumerate}
\item Monodromy maps for open book decompositions supporting tight contact structures which are not holomorphically fillable.
\item Right-veering monodromy maps of open book decompositions supporting overtwisted contact structures.

\item Explicit examples  on the punctured torus.

\end{enumerate}

We briefly explain each:  

(1) There are tight contact structures which are
not weakly symplectically fillable and also weakly symplectically fillable ones
which are not strongly fillable (or holomorphically fillable).  Contact structures
which are tight but are not fillable were first discovered by Etnyre-Honda~\cite{EH}
and expanded to infinitely many examples in a series of papers by Lisca-Stipsicz
(see \cite{LS1,LS2} for the first couple).
Weakly fillable contact structures which are not strongly fillable were first exhibited
by Eliashberg~\cite{El2}; further examples were given by Ding-Geiges~\cite{DG}.  
Let $(S,h)$ be any open book decomposition for a contact structure $(M,\xi)$ which
is tight but not holomorphically fillable.  Then $h$ must be right-veering
by Theorem~\ref{veer} but cannot be a product of positive Dehn twists
by Corollary~\ref{holo}.    

(2) As explained in Proposition~\ref{rvstabilization},
any overtwisted contact structure admits an open book $(S,h)$ which
is right-veering.  Such a monodromy map $h$ can never be a
product of positive Dehn twists, again by Corollary~\ref{holo}.   

(3) is the topic of the sequel~\cite{HKM2}.

\s\n
{\em Acknowledgements.} We thank Francis Bonahon and Bill Thurston for helpful discussions.

\end{document}